\documentclass[twoside,12pt,draft]{article}
\usepackage{amssymb}

\usepackage{makeidx}
\usepackage{a4wide}
\usepackage{amsthm}
\usepackage{amsmath}
\usepackage[all]{xy}
\usepackage{enumerate}
\usepackage{fancyhdr}

\newtheorem{theorem}{Theorem}[section]
\newtheorem{proposition}[theorem]{Proposition}
\newtheorem{corollary}[theorem]{Corollary}
\newtheorem{lemma}[theorem]{Lemma}
\theoremstyle{definition}
\newtheorem{example}[theorem]{Example}
\newtheorem{c-example}[theorem]{Counter Example}
\newtheorem*{Notation}{Notation}
\newtheorem*{Beweis}{Proof}
\newtheorem{definition}[theorem]{Definition}
\newtheorem{punto}[theorem]{}
\theoremstyle{remark}
\newtheorem{remark}[theorem]{Remark}

\CompileMatrices
\input{tcilatex}

\begin{document}

\title{Dual Entwining Structures and Dual Entwined Modules \thanks{%
MSC 2000: 16W30, 18E15 \newline
Keywords: Entwining structures, Entwined modules, Doi-Koppinen structures -
Doi-Koppinen modules, Hopf-Galois (co)extensions, duality, (co)module
(co)algebras}}
\author{\textbf{Jawad Y. Abuhlail} \\
Mathematics Department\\
Birzeit University\\
P.O. Box 14, Birzeit - Palestine \\
jabuhlail@birzeit.edu}
\date{}
\maketitle

\begin{abstract}
In this note we introduce and investigate the concepts of dual entwining
structures and dual entwined modules. This generalizes the concepts of dual
Doi-Koppinen structures and dual Doi-Koppinen modules introduced (in the
infinite case over rings) by the author is his dissertation.
\end{abstract}

\section*{Introduction}

This note deals with the following problem: let $(A,C,\psi )$ be a given
entwining structure over a commutative base ring $R.$ Find an $R$-subalgebra 
$\widetilde{A}\subseteq C^{\ast }$ and an $R$-coalgebra $\widetilde{C}%
\subseteq A^{\ast },$ such that $(\widetilde{A},\widetilde{C},\psi ^{\ast })$
is an entwining structure.

For general entwining structures, it's not clear if such a \emph{dual }%
entwining structure exists. However, once it's found, we have the \emph{%
expected} duality theorems between the corresponding categories of entwined
modules. For the special case of Doi-Koppinen structures over noetherian
rings, the problem was solved by the author in his dissertation. Our results
are formulated for right-right entwining structures. Corresponding versions
for left-left, right-left and left-right entwining structures can be derived
easily using the left-right dictionary (e.g. \cite{CMZ2002}).

The paper consists of three sections. In the first section, we give the
necessary definitions and results from the theory of Hopf algebras and
entwining structures. In the second section we present and investigate the
concepts of dual entwining structures and dual entwined modules. The third
sections is an extended version of \cite[Paragraph 3.4]{Abu2001} formulated
for \emph{right-right }Doi-Koppinen structures. Our results in the third
section generalize also those \emph{achieved independently} by L. Zhang \cite
{Zha98} on \emph{dual relative Hopf modules} in the case of a commutative
base field.

Throughout this paper $R$ denotes a commutative ring with $1_{R}\neq 0_{R}.$
The category of $R$-(bi)modules will be denoted by $\mathcal{M}_{R}.$ For an 
$R$-coalgebra $(C,\Delta _{C},\varepsilon _{C})$ and an $R$-algebra $(A,\mu
_{A},\eta _{A})$ we consider $\mathrm{Hom}_{R}(C,A)\;$as an $R$-algebra with
the so called \emph{convolution product }$(f\star g)(c):=\sum
f(c_{1})g(c_{2})$ and unity $\eta _{A}\circ \varepsilon _{C}.$ For an $R$%
-algebra $A$ and an $A$-module $M,$ an $A$-submodule $N\subset M$ will be
called $R$\emph{-cofinite}, if $M/N$ is f.g. in $\mathcal{M}_{R}.$ We call
an $R$-submodule $K\subseteq M$ \emph{pure} (in the sense of Cohn), if the
canonical map $\iota _{K}\otimes id_{N}:K\otimes _{R}N\rightarrow M\otimes
_{R}N$ is injective for every $R$-module $N.$

\section{Preliminaries}

In this section we give some definitions and Lemmata from the theory of Hopf
Algebras and Entwining Structures.

\begin{punto}
\textbf{Measuring }$R$-\textbf{pairings. }If $C$ is an $R$-coalgebra and $A$
is an $R$-algebra with a morphism of $R$-algebras $\kappa :A\rightarrow $ $%
C^{\ast },$ $a\mapsto \lbrack c\mapsto <a,c>],$ then we call $P:=(A,C)$ a 
\emph{measuring }$R$\emph{-pairing}. In this case $C$ is an $A$-bimodule
through 
\begin{equation}
a\rightharpoonup c:=\sum c_{1}<a,c_{2}>\text{ and }c\leftharpoonup a:=\sum
<a,c_{1}>c_{2}\text{ for all }a\in A,\text{ }c\in C.  \label{C-r}
\end{equation}
Let $(A,C)$ and $(B,D)$ be measuring $R$-pairings, $\xi :A\rightarrow B$ an $%
R$-algebra morphism and $\theta :D\rightarrow C$ an $R$-coalgebra morphism.
Then we say $(\xi ,\theta ):(B,D)\rightarrow (A,C)$ is a \emph{morphism of
measuring }$R$\emph{-pairings,} if 
\begin{equation*}
<\xi (a),d>=<a,\theta (d)>\text{ for all }a\in A\text{ and }d\in D.
\end{equation*}
\newline
The category of measuring $R$-pairings and morphisms described above will be
denoted by $\mathcal{P}_{m}.$
\end{punto}

\begin{punto}
\textbf{The }$\alpha $\textbf{-condition}. Let $P=(A,C)$ be a measuring $R$%
-pairing. We say $P$ satisfies the $\alpha $\emph{-condition}\textbf{\ }(or $%
P$ is a \emph{measuring }$\alpha $\emph{-pairing}), if for every $R$-module $%
M$ the following map is injective: 
\begin{equation}
\alpha _{M}^{P}:M\otimes _{R}C\ \rightarrow \mathrm{Hom}_{R}(A,M),\text{ }%
\sum m_{i}\otimes c_{i}\mapsto \lbrack a\mapsto \sum m_{i}<a,c_{i}>].
\label{alp}
\end{equation}
With $\mathcal{P}_{m}^{\alpha }\subset \mathcal{P}_{m}$ we denote the \emph{%
full }subcategory of measuring $\alpha $-pairings.

We say an $R$-coalgebra $C$ satisfies the $\alpha $\emph{-condition,} if $%
(C^{\ast },C)$ satisfies the $\alpha $-condition (equivalently, if $_{R}C$
is \emph{locally projective} in the sense of B. Zimmermann-Huignes 
\cite[Theorem 2.1]{Z-H76}, \cite[Theorem 3.2]{Gar76}).
\end{punto}

\begin{punto}
\textbf{Subgenerators}. Let $A$ be an $R$-algebra and $K$ an $A$-module. We
say an $A$-module $N$ is $K$\emph{-subgenerated}, if $N$ is isomorphic to a
submodule of a $K$-generated $A$-module (equivalently, if $N$ is Kernel of $K
$-generated $A$-modules). The \emph{full }subcategory of $A$-modules, whose
objects are the $K$-subgenerated $A$-modules is denoted by $\sigma \lbrack
K].$ Moreover $\sigma \lbrack K]$ is the \emph{smallest }Grothendieck full
subcategory of the category of $A$-modules that contains $K.$ The reader is
referred to \cite{Wis88} for the well developed theory of categories of this
type.
\end{punto}

\section*{\textbf{Rational Modules}}

\begin{punto}
\label{rat-dar}Let $P=(A,C)$ a measuring $\alpha $-pairing. Let $M$ be a
left (a right) $A$-module, $\rho _{M}:M\rightarrow \mathrm{Hom}_{R}(A,M)$
the canonical $A$-linear map and put $\mathrm{Rat}^{C}(_{A}M):=\rho
_{M}^{-1}(M\otimes _{R}C)$ (resp. $^{C}\mathrm{Rat}(M_{A}):=\rho
_{M}^{-1}(C\otimes _{R}M)$). We call $_{A}M$ (resp. $M_{A}$) $C$\emph{%
-rational,} if $\mathrm{Rat}^{C}(_{A}M)=M$ (resp. $^{C}\mathrm{Rat}(M_{A})=M$%
). If $M$ is an $A$-bimodule, then we set $^{C}\mathrm{Rat}^{C}(_{A}M_{A})=%
\mathrm{Rat}^{C}(_{A}M)\cap $ $^{C}\mathrm{Rat}(M_{A})$ and call $M$ $C$%
-birational, if $^{C}\mathrm{Rat}^{C}(_{A}M_{A})=M.$
\end{punto}

\begin{lemma}
\label{clos}\emph{(\cite[Lemma 2.2.7]{Abu2001})}\ Let $P=(A,C)$ be a
measuring $\alpha $-pairing. For every left (resp. right) $A$-module $M$ we
have:

\begin{enumerate}
\item  $\mathrm{Rat}^{C}(_{A}M)\subset M$ \emph{(}resp. $^{C}\mathrm{Rat}%
(M_{A})\subset M$\emph{)} is an $A$-submodule.

\item  For every $A$-submodule $N\subset M,$ it follows that $\mathrm{Rat}%
^{C}(_{A}N)=N\cap \mathrm{Rat}^{C}(_{A}M)$ \emph{(}resp. $^{C}\mathrm{Rat}%
(N_{A})=N\cap $ $^{C}\mathrm{Rat}(M_{A})$\emph{)}.

\item  $\mathrm{Rat}^{C}(\mathrm{Rat}^{C}(_{A}M))=\mathrm{Rat}^{C}(_{A}M)$ 
\emph{(}resp. $^{C}\mathrm{Rat}(^{C}\mathrm{Rat}(M_{A}))=$ $^{C}\mathrm{Rat}%
(M_{A})$\emph{)}.

\item  For a \emph{(}resp. a right\emph{)} $A$-module $L$ and an $A$-linear
map $f:$ $M\rightarrow L,$ we have $f(\mathrm{Rat}^{C}(_{A}M))\subseteq 
\mathrm{Rat}^{C}(_{A}L)$ \emph{(}resp. $f(^{C}\mathrm{Rat}(M_{A}))\subseteq $
$^{C}\mathrm{Rat}(L_{A})$\emph{)}.
\end{enumerate}
\end{lemma}

\begin{Notation}
For a measuring $\alpha $-pairing $(A,C)$ we denote with $\mathrm{Rat}%
^{C}(_{A}\mathcal{M})\subseteq $ $_{A}\mathcal{M}$ (resp. $^{C}\mathrm{Rat}(%
\mathcal{M}_{A})\subseteq \mathcal{M}_{A},$ $^{C}\mathrm{Rat}^{C}(_{A}%
\mathcal{M}_{A})\subseteq $ $_{A}\mathcal{M}_{A}$) the \emph{full }%
subcategory of $C$-rational left $A$-modules (resp. $C$-rational right $A$%
-modules, $C$-birational $A$-bimodules).
\end{Notation}

\begin{theorem}
\label{equal}\emph{(\cite[Satz 2.2.16, Folgerung 2.2.22]{Abu2001})} For a
measuring $R$-pairing $P=(${$A$}$,${$C$}$)$ the following are equivalent:

\begin{enumerate}
\item  $P$ satisfies the $\alpha $-condition;

\item  $_{R}C$ is locally projective and $\kappa _{P}(A)\subseteq C^{\ast }$
is dense \emph{(}w.r.t. the finite topology\emph{).}

If these equivalent conditions are satisfied, then we have isomorphisms of
categories 
\begin{equation*}
\begin{tabular}{lllll}
$\mathcal{M}^{{C}}$ & $\simeq $ & $\mathrm{Rat}^{{C}}(_{A}\mathcal{M})$ & $=$
& $\sigma \lbrack _{{A}}{C}]$ \\ 
& $\simeq $ & $\mathrm{Rat}^{{C}}(_{C^{\ast }}\mathcal{M})$ & $=$ & $\sigma
\lbrack _{C^{\ast }}{C}];$ \\ 
$^{C}\mathcal{M}$ & $\simeq $ & $^{{C}}\mathrm{Rat}(\mathcal{M}_{A})$ & $=$
& $\sigma \lbrack {C}_{{A}}]$ \\ 
& $\simeq $ & $^{{C}}\mathrm{Rat}(\mathcal{M}_{C^{\ast }})$ & $=$ & $\sigma
\lbrack {C}_{C^{\ast }}];$ \\ 
$^{C}\mathcal{M}^{C}$ & $\simeq $ & $^{{C}}\mathrm{Rat}^{{C}}(_{A}\mathcal{M}%
_{A})$ & $=$ & $\sigma \lbrack _{A}({C\otimes _{R}C)}_{{A}}]$ \\ 
& $\simeq $ & $^{{C}}\mathrm{Rat}^{{C}}(_{C^{\ast }}\mathcal{M}_{C^{\ast }})$
& $=$ & $\sigma \lbrack _{C^{\ast }}({C\otimes _{R}C)}_{C^{\ast }}]$%
\end{tabular}
\end{equation*}
\end{enumerate}
\end{theorem}

\begin{lemma}
\label{prop-mes}\emph{(\cite[Lemma 2.1.23]{Abu2001}) }Let $P=(A,C),$ $%
Q=(B,D) $ be measuring $R$-pairings and $\xi :A\rightarrow B,$ $\theta
:D\rightarrow C$ be $R$-linear maps with 
\begin{equation*}
<\xi (a),d>=<a,\theta (d)>\text{ for all }a\in A\text{ and }d\in D.
\end{equation*}

\begin{enumerate}
\item  Set $P\otimes P:=(A\otimes _{R}A,C\otimes _{R}C)$ and assume that $%
C\otimes _{R}C\overset{\chi _{P\otimes P}}{\hookrightarrow }(A\otimes
_{R}A)^{\ast }$ is an embedding. If $\xi $ is an $R$-algebra morphism, then $%
\theta $ is an $R$-coalgebra morphism. Moreover, if $A$ is commutative, then 
$C$ is cocommutative.

\item  Assume $B\overset{\kappa _{Q}}{\hookrightarrow }D^{\ast }$ to be an
embedding. If $\theta $ is an $R$-coalgebra morphism, then $\xi $ is an $R$%
-algebra morphism. Moreover, if $C$ is cocommutative and $A\subseteq C^{\ast
},$ then $A$ is commutative.
\end{enumerate}
\end{lemma}

\begin{punto}
\label{du-co}(\cite[Theorem 2.8]{AG-TW2000}, \cite[Remark 2.14, Prosposition
2.15]{AG-TL2001}) Assume $R$ to be noetherian. Let $A$ be an $R$-algebra and
consider $A^{\ast }$ as an $A$-bimodule through the left and the right
regular $A$-action $(af)(b)=f(ba)$ and $(fa)(b)=f(ab).$ We define the \emph{%
finite dual} of $A$ as the $R$-module 
\begin{equation*}
\begin{tabular}{lll}
$A^{\circ }$ & $:=$ & $\{f\in A^{\ast }\mid AfA\text{ is f.g. in }\mathcal{M}%
_{R}\}$ \\ 
& $=$ & $\{f\in A^{\ast }\mid $ $f(I)=0$ for some $R$-cofinite ideal $%
I\vartriangleleft A\}.$%
\end{tabular}
\end{equation*}
An $R$-algebra (resp. an $R$-bialgebra, a Hopf $R$-algebra) $A$ with $%
A^{\circ }\subset R^{A}$ pure will be called an $\alpha $\emph{-algebra} (an 
$\alpha $\emph{-bialgebra}, a \emph{Hopf }$\alpha $\emph{-algebra}). For
every $\alpha $-algebra (resp. $\alpha $-bialgebra, Hopf $\alpha $-algebra) $%
A,$ the finite dual $A^{\circ }$ becomes a locally projective $R$-coalgebra
(resp. $R$-bialgebra, Hopf $R$-algebra). If $A$ is an $\alpha $-algebra and $%
\widetilde{C}\subseteq A^{\circ }$ is an $R$-subcoalgebra, then $(A,%
\widetilde{C})$ is a measuring $\alpha $-pairing. For $\alpha $-algebras
(resp. $\alpha $-bialgebra, Hopf $\alpha $-algebras) $A,B$ and a morphism of 
$R$-algebras (resp. $R$-bialgebras, Hopf $R$-algebras) $f:A\rightarrow B,$
it follows directly from Lemma \ref{prop-mes} that the restriction of $%
f^{\ast }:B^{\ast }\rightarrow A^{\ast }$ to $B^{\circ }$ induces a morphism
of $R$-coalgebras (resp. $R$-bialgebras, Hopf $R$-algebras) $f^{\circ
}:B^{\circ }\rightarrow A^{\circ }.$
\end{punto}

\begin{remark}
\label{emb}(\cite[Folgerung 2.1.10 (1)]{Abu2001})\emph{\ }Let $V,W$ be $R$%
-modules and $X\subseteq V^{\ast },$ $Y\subseteq W^{\ast }$ be $R$%
-submodules. If $R$ is noetherian and $X\subset R^{V}$ is $Y$-pure (or $%
Y\subset R^{W}$\ is $X$-pure), then the following induced canonical map is
injective: 
\begin{equation*}
\varpi :X\otimes _{R}Y\rightarrow (V\otimes _{R}W)^{\ast },\text{ }f\otimes
g\mapsto f\underline{\otimes }g,\text{ where }(f\underline{\otimes }%
g)(v\otimes w):=f(v)g(w).
\end{equation*}
\end{remark}

\section*{\textbf{Entwined modules}}

\begin{punto}
\label{ent-str}\ A \emph{right-right entwining structure}\textbf{\ }$%
(A,C,\psi )$ consists of an $R$-algebra $(A,\mu _{A},\eta _{A}),$ an $R$%
-coalgebra $(C,\Delta _{C},\varepsilon _{C})$ and an $R$-linear map 
\begin{equation*}
\psi :C\otimes _{R}A\rightarrow A\otimes _{R}C,\text{ }c\otimes a\mapsto
\sum a_{\psi }\otimes c^{\psi },
\end{equation*}
such that 
\begin{equation}
\begin{tabular}{llllll}
$\sum (ab)_{\psi }\otimes c^{\psi }$ & $=$ & $\sum a_{\psi }b_{\Psi }\otimes
c^{\psi \Psi },$ & $\sum (1_{A})_{\psi }\otimes c^{\psi }$ & $=$ & $%
1_{A}\otimes c,$ \\ 
$\sum a_{\psi }\otimes \Delta _{C}(c^{\psi })$ & $=$ & $\sum a_{\psi \Psi
}\otimes c_{1}^{\Psi }\otimes c_{2}^{\psi },$ & $\sum \varepsilon
_{C}(c^{\psi })a_{\psi }$ & $=$ & $\varepsilon _{C}(c)a.$%
\end{tabular}
\label{rr-ent}
\end{equation}
Let $(A,C,\psi )$ and $(B,D,\Psi )$ be right-right entwining structures. A
morphism $(\gamma ,\delta ):(A,C,\psi )\rightarrow (B,D,\Psi )$ consists of
an $R$-algebra morphism $\gamma :A\rightarrow B$ and an $R$-coalgebra
morphism $\delta :C\rightarrow D,$ such that 
\begin{equation*}
\sum \gamma (a_{\psi })\otimes \delta (c^{\psi })=\sum \gamma (a)_{\Psi
}\otimes \delta (c)^{\Psi }.
\end{equation*}
With $\mathbb{E}_{\bullet }^{\bullet }$ we denote the category of
right-right entwining structures. For definitions of the categories of
left-left, right-left and left-right entwining structures the interested
reader may refer to \cite{CMZ2002}.
\end{punto}

\begin{punto}
Let $(A,C,\psi )$ be a right-right entwining structure. An \emph{entwined
module} corresponding to $(A,C,\psi )$ is a right $A$-module $M,$ which is
also a right $C$-comodule, such that 
\begin{equation*}
\varrho _{M}(ma)=\sum m_{<0>}a_{\psi }\otimes m_{<1>}^{\psi }\text{ for all }%
m\in M,\text{ }a\in A.
\end{equation*}
For entwined modules $M,N$ corresponding to $(A,C,\psi )$ we denote with $%
\mathrm{Hom}_{A}^{C}(M,N)$ the set of $A$-linear $C$-colinear morphisms from 
$M$ to $N.$ The category of right-right entwined modules and $A$-linear $C$%
-colinear morphisms is denoted by $\mathcal{M}_{A}^{C}(\psi ).$ Entwined
modules were introduced by T. Brzezi\'{n}ski and S. Majid in \cite{BM98} as
a generalization of Doi-Koppinen modules presented in \cite{Doi92} and \cite
{Kop95}.
\end{punto}

\begin{lemma}
\label{lem-entw}Let $(A,C,\psi )$ be a right-right entwining structure over $%
R$ and set $\mathcal{C}:=A\otimes _{R}C.$

\begin{enumerate}
\item  $\mathcal{C}$ is an $A$-coring with $A$-bimodules structure given by 
\begin{equation}
a(\widetilde{a}\otimes c):=a\widetilde{a}\otimes c,\text{ }(\widetilde{a}%
\otimes c)a:=\sum \widetilde{a}a_{\psi }\otimes c^{\psi },  \label{AotC-mod}
\end{equation}
comultiplication 
\begin{equation*}
\Delta _{\mathcal{C}}:A\otimes _{R}C\rightarrow (A\otimes _{R}C)\otimes
_{A}(A\otimes _{R}C),\text{ }a\otimes c\mapsto \sum (a\otimes c_{1})\otimes
_{A}(1_{A}\otimes c_{2})
\end{equation*}
and counity $\varepsilon _{\mathcal{C}}:=\vartheta _{A}^{r}\circ
(id_{A}\otimes \varepsilon _{C}).$ Moreover $\mathcal{M}_{A}^{C}(\psi
)\simeq \mathcal{M}^{\mathcal{C}}.$

\item  $\#_{\psi }^{op}(C,A):=\mathrm{Hom}_{R}(C,A)$ is an $A$-ring with $A$%
-bimodule structure given by $(af)(c):=\sum a_{\psi }f(c^{\psi }),$ $%
(fa)(c):=f(c)a,$ multiplication 
\begin{equation}
(f\cdot g)(c)=\sum f(c_{2})_{\psi }g(c_{1}^{\psi })  \label{re-m}
\end{equation}
and unity $\eta _{A}\circ \varepsilon _{C}.$

\item  Consider $^{\ast }\mathcal{C}:=\mathrm{Hom}_{A-}(\mathcal{C},A)$ as
an $A$-ring with the canonical $A$-bimodule structure, multiplication 
\begin{equation*}
(f\star _{l}g)(c)=\sum g(c_{1}f(c_{2}))\text{ for all }f,g\in \text{ }^{\ast
}\mathcal{C}\text{ and }c\in \mathcal{C}
\end{equation*}
and unity $\varepsilon _{\mathcal{C}}.$ Then $\#_{\psi }^{op}(C,A)\simeq $ $%
^{\ast }\mathcal{C}$ as $A$-rings via 
\begin{equation}
\nu :\mathrm{Hom}_{R}(C,A)\longrightarrow \mathrm{Hom}_{A-}(A\otimes
_{R}C,A),\text{ }f\mapsto \lbrack a\otimes c\mapsto af(c)]  \label{vphi-iso}
\end{equation}
with inverse $h\mapsto \lbrack c\mapsto h(1_{A}\otimes c)].$
\end{enumerate}
\end{lemma}

\begin{Beweis}
\begin{enumerate}
\item  This was noticed first by M. Takeuchi and can be found in several
references (e.g. \cite[Proposition 2.2.]{Brz}).

\item  For all $a,b\in A,$ $f\in \#_{\psi }^{op}(C,A)$ and $c\in C$ we have 
\begin{equation*}
\begin{tabular}{lllllll}
$((ab)f)(c)$ & $=$ & $\sum (ab)_{\psi }f(c^{\psi })$ & $=$ & $\sum a_{\psi
}b_{\Psi }f(c^{\psi \Psi })$ &  &  \\ 
& $=$ & $\sum a_{\psi }(bf)(c^{\psi })$ & $=$ & $(a(bf))(c).$ &  & 
\end{tabular}
\end{equation*}

It's clear then that the left and the right $A$-actions given above define
on $\#_{\psi }^{op}(C,A)$ a structure of an $A$-bimodule. Moreover we have
for all $f,g,h\in \mathrm{Hom}_{R}(C,A)$ and $c\in C:$%
\begin{equation*}
\begin{tabular}{llll}
$((f\cdot g)\cdot h)(c)$ & $=$ & $\sum [(f\cdot g)(c_{2})]_{\widehat{\psi }%
}h(c_{1}^{\widehat{\psi }})$ &  \\ 
& $=$ & $\sum [f(c_{22})_{\psi }g(c_{21}^{\psi })]_{\widehat{\psi }}h(c_{1}^{%
\widehat{\psi }})$ &  \\ 
& $=$ & $\sum [f(c_{22})_{\psi \widehat{\psi }}g(c_{21}^{\psi })_{\psi
^{\prime }}]h((c_{1}^{\widehat{\psi }})^{\psi ^{\prime }})$ &  \\ 
& $=$ & $\sum [f(c_{2})_{\psi \widehat{\psi }}g(c_{12}^{\psi })_{\psi
^{\prime }}]h((c_{11}^{\widehat{\psi }})^{\psi ^{\prime }})$ &  \\ 
& $=$ & $\sum f(c_{2})_{\psi }g((c_{1}^{\psi })_{2})_{\psi ^{\prime
}}h((c_{1}^{\psi })_{1}^{\psi ^{\prime }})$ &  \\ 
& $=$ & $\sum f(c_{2})_{\psi }(g\cdot h)(c_{1}^{\psi })$ &  \\ 
& $=$ & $(f\cdot (g\cdot h))(c).$ & 
\end{tabular}
\end{equation*}

It's clear that $\eta _{A}\circ \varepsilon _{C}$ is a unity for $\#_{\psi
}^{op}(C,A).$

\item  Note that $\nu $ is given by the canonical isomorphisms 
\begin{equation*}
\mathrm{Hom}_{R}(C,A)\simeq \mathrm{Hom}_{R}(C,\mathrm{Hom}_{A-}(A,A))\simeq 
\mathrm{Hom}_{A-}(A\otimes _{R}C,A).
\end{equation*}
For all $a\in A,$ $f\in \#_{\psi }^{op}(C,A)\;$and $c\in C$ we have 
\begin{equation*}
\begin{tabular}{lllll}
$\nu (af)(b\otimes c)$ & $=$ & $b((af)(c))$ & $=$ & $b(\sum a_{\psi
}f(c^{\psi }))$ \\ 
& $=$ & $\sum ba_{\psi }f(c^{\psi })$ & $=$ & $\nu (f)(\sum ba_{\psi
}\otimes c^{\psi })$ \\ 
& $=$ & $\nu (f)((b\otimes c)a)$ & $=$ & $(a\nu (f))(b\otimes c).$%
\end{tabular}
\end{equation*}
It is obvious that $\nu $ is right $A$-linear. For all $f,g\in \#_{\psi
}^{op}(C,A),$ $a\in A$ and $c\in C$ we have 
\begin{equation*}
\begin{tabular}{lllll}
$\nu (f\cdot g)(a\otimes c)$ & $=$ & $a((f\cdot g)(c))$ & $=$ & $a\sum
f(c_{2})_{\psi }g(c_{1}^{\psi })$ \\ 
& $=$ & $\sum af(c_{2})_{\psi }g(c_{1}^{\psi })$ & $=$ & $\nu (g)(\sum
af(c_{2})_{\psi }\otimes c_{1}^{\psi })$ \\ 
& $=$ & $\nu (g)(\sum (a\otimes c_{1})f(c_{2}))$ & $=$ & $\nu (g)(\sum
(a\otimes c_{1})1_{A}f(c_{2}))$ \\ 
& $=$ & $\nu (g)(\sum (a\otimes c_{1})\nu (f)(1_{A}\otimes c_{2}))$ & $=$ & $%
(\nu (f)\star _{l}\nu (g))(a\otimes c).$%
\end{tabular}
\end{equation*}
Consequently $\nu $ is an isomorphism of $A$-rings.$\blacksquare $
\end{enumerate}
\end{Beweis}

\begin{punto}
\label{s-rat}Let $(A,C,\psi )$ be a right-right entwining structure over $R$
and consider the corresponding $A$-coring $\mathcal{C}:=A\otimes _{R}C.$ We
say that $(A,C,\psi )$ \emph{satisfies the }$\alpha $\emph{-condition} (or
is an $\alpha $\emph{-entwining structure}) if for every right $A$-module $%
M, $ the following map is injective 
\begin{equation*}
\alpha _{M}^{\psi }:M\otimes _{R}C\rightarrow \mathrm{Hom}_{R}(\#_{\psi
}^{op}(C,A),M),\text{ }m\otimes c\mapsto \lbrack f\mapsto mf(c)]
\end{equation*}
(equivalently if $_{A}\mathcal{C}$ is locally projective).
\end{punto}

Inspired by \cite[3.1]{Doi94} we introduce

\begin{definition}
Let $(A,C,\psi )$ be a right-right entwining structure that satisfies the $%
\alpha $-condition. Let $M$ be a right $\#_{\psi }^{op}(C,A)$-module, $\rho
_{M}:M\rightarrow \mathrm{Hom}_{-\#_{\psi }^{op}(C,A)}(\#_{\psi
}^{op}(C,A),M)$ the canonical map and put $\mathrm{Rat}^{C}(M_{\#_{\psi
}^{op}(C,A)}):=\rho _{M}^{-1}(M\otimes _{R}C).$ Then $M$ will be called $\#$%
\emph{-rational}, if $\mathrm{Rat}^{C}(M_{\#_{\psi }^{op}(C,A)})=M.$ For a $%
\#$-rational right $\#_{\psi }^{op}(C,A)$-module $M$ we set $\varrho
_{M}:=(\alpha _{M}^{\psi })^{-1}\circ \rho _{M}:M\rightarrow M\otimes _{R}C.$
The category of $\#$-rational right $\#_{\psi }^{op}(C,A)$-modules and $%
\#_{\psi }^{op}(C,A)$-linear maps will be denoted with $\mathrm{Rat}^{C}(%
\mathcal{M}_{\#_{\psi }^{op}(C,A)}).$
\end{definition}

\begin{theorem}
\label{ent-sg}\emph{(\cite[Lemma 3.8, Theorem 3.10]{Abu})}\ Let $(A,C,\psi )$
be a right-right entwining structure and consider the corresponding $A$%
-coring $\mathcal{C}:=A\otimes _{R}C.$

\begin{enumerate}
\item  If $_{R}C$ is flat, then $_{A}\mathcal{C}$ is flat and $\mathcal{M}%
_{A}^{C}(\psi )$ is a Grothendieck category with enough injective objects.

\item  If $_{R}C$ is locally projective \emph{(}resp. f.g. projective\emph{)}%
, then $_{A}\mathcal{C}$ is locally projective \emph{(}resp.\emph{\ }f.g.
projective\emph{)}\ and 
\begin{equation}
\mathcal{M}_{A}^{C}(\psi )\simeq \mathrm{Rat}^{C}(\mathcal{M}_{\#_{\psi
}^{op}(C,A)})\simeq \sigma \lbrack (A\otimes _{R}C)_{\#_{\psi }^{op}(C,A)}]%
\text{\ \emph{(}resp. }\mathcal{M}_{A}^{C}(\psi )\simeq \mathcal{M}%
_{\#_{\psi }^{op}(C,A)}\text{\emph{)}.}  \label{iso-sm}
\end{equation}
\end{enumerate}
\end{theorem}

\section{\textbf{Dual entwined modules}}

In this section we fix the following:\ $R$ is noetherian, $(A,C,\psi )$ is a
right-right entwining structure with $A$ an $\alpha $-algebra and $%
\widetilde{A}\subseteq C^{\ast }$ is an $R$-subalgebra with $\varepsilon
_{C}\in \widetilde{A},$ $\widetilde{C}\subseteq A^{\circ }$ is an $R$%
-subcoalgebra. So we have a measuring $R$-pairing $(\widetilde{A},C)$ and a
measuring $\alpha $-pairing $(A,\widetilde{C}).$ Besides the above technical
assumptions we assume moreover that $\psi ^{\ast }(\widetilde{C}\otimes _{R}%
\widetilde{A})\subseteq \widetilde{A}\otimes _{R}\widetilde{C},$ i.e. the
following diagram 
\begin{equation*}
\xymatrix{ (A \otimes_R C)^* \ar[rr]^{\psi ^*} & & (C \otimes_R A)^* \\
{\widetilde{C}} \otimes_{R} {\widetilde{A}} \ar@{.>}[rr]_{\varphi}
\ar@{^{(}->}[u]& & {\widetilde{A}} \otimes_{R} {\widetilde{C}}
\ar@{^{(}->}[u]}
\end{equation*}
can be completed commutatively with an $R$-linear morphism 
\begin{equation}
\varphi :\widetilde{C}\otimes _{R}\widetilde{A}\rightarrow \widetilde{A}%
\otimes _{R}\widetilde{C},\text{ }\widetilde{f}\otimes \widetilde{g}\mapsto
\sum \widetilde{g}_{\varphi }\underline{\otimes }\widetilde{f}^{\varphi },%
\text{ where }(\sum \widetilde{g}_{\varphi }\underline{\otimes }\widetilde{f}%
^{\varphi })(c\otimes a):=\sum \widetilde{f}(a_{\psi })\widetilde{g}(c^{\psi
}).  \label{varphi}
\end{equation}

\begin{theorem}
$(\widetilde{A},\widetilde{C},\varphi )$ is a right-right entwining
structure and we have isomorphisms of categories 
\begin{equation}
\mathcal{M}_{\widetilde{A}}^{\widetilde{C}}(\varphi )\simeq \mathrm{Rat}^{%
\widetilde{C}}(\mathcal{M}_{\#_{\varphi }^{op}(\widetilde{C},\widetilde{A}%
)})=\sigma \lbrack (\widetilde{A}\otimes _{R}\widetilde{C})_{\#_{\varphi
}^{op}(\widetilde{C},\widetilde{A})}].  \label{du-ent-is}
\end{equation}
If moreover $_{R}\widetilde{C}$ is f.g. projective, then 
\begin{equation}
\mathcal{M}_{\widetilde{A}}^{\widetilde{C}}(\varphi )\simeq \mathcal{M}%
_{\#_{\varphi }^{op}(\widetilde{C},\widetilde{A})}).  \label{du-ent-is-fg}
\end{equation}
\end{theorem}

\begin{Beweis}
Let $\widetilde{f}\in \widetilde{C},$ $\widetilde{g},$ $\widetilde{h}\in 
\widetilde{A},$ $c\in C$ and $a,b\in A$ be arbitrary. Then we have 
\begin{equation*}
\begin{tabular}{lll}
$\sum ((\widetilde{g}\star \widetilde{h})_{\varphi }\underline{\otimes }%
\widetilde{f}^{\varphi })(c\otimes a)$ & $=$ & $\sum \widetilde{f}(a_{\psi
})(\widetilde{g}\star \widetilde{h})(c^{\psi })$ \\ 
& $=$ & $\sum \widetilde{f}(a_{\psi })\widetilde{g}((c^{\psi })_{1})%
\widetilde{h}((c^{\psi })_{2})$ \\ 
& $=$ & $\sum \widetilde{f}(a_{\psi \Psi })\widetilde{g}(c_{1}^{\Psi })%
\widetilde{h}(c_{2}^{\psi })$ \\ 
& $=$ & $\sum (\widetilde{g}_{\varphi }\underline{\otimes }\widetilde{h}%
_{\Phi }\underline{\otimes }\widetilde{f}^{\varphi \Phi })(c_{1}\otimes
c_{2}\otimes a)$ \\ 
& $=$ & $\sum ((\widetilde{g}_{\varphi }\star \widetilde{h}_{\Phi })%
\underline{\otimes }\widetilde{f}^{\varphi \Phi })(c\otimes a)$%
\end{tabular}
\end{equation*}
and 
\begin{equation*}
(\sum (\varepsilon _{C})_{\varphi }\underline{\otimes }\widetilde{f}%
^{\varphi })(c\otimes a)=\sum \widetilde{f}(a_{\psi })\varepsilon
_{C}(c^{\psi })=\widetilde{f}(\varepsilon _{C}(c)a)=(1_{\widetilde{A}}%
\underline{\otimes }\widetilde{f})(c\otimes a).
\end{equation*}
On the other hand we have 
\begin{equation*}
\begin{tabular}{lll}
$(\sum \widetilde{g}_{\varphi }\underline{\otimes }((\widetilde{f}^{\varphi
})_{1}\underline{\otimes }(\widetilde{f}^{\varphi })_{2}))(c\otimes a\otimes
b)$ & $=$ & $(\sum \widetilde{g}_{\varphi }\underline{\otimes }\widetilde{f}%
^{\varphi })(c\otimes ab)$ \\ 
& $=$ & $\sum \widetilde{f}((ab)_{\psi })\widetilde{g}(c^{\psi })$ \\ 
& $=$ & $\sum \widetilde{f}(a_{\psi }b_{\Psi })\widetilde{g}(c^{\psi \Psi })$
\\ 
& $=$ & $\sum \widetilde{f}_{1}(a_{\psi })\widetilde{f}_{2}(b_{\Psi })%
\widetilde{g}(c^{\psi \Psi })$ \\ 
& $=$ & $\sum (\widetilde{g}_{\varphi \Phi }\underline{\otimes }\widetilde{f}%
_{1}^{\Phi }\underline{\otimes }\widetilde{f}_{2}^{\varphi })(c\otimes
a\otimes b)$%
\end{tabular}
\end{equation*}
and 
\begin{equation*}
(\sum \varepsilon _{\widetilde{C}}(\widetilde{f}^{\varphi })\widetilde{g}%
_{\varphi })(c)=\sum (\widetilde{g}_{\varphi }\underline{\otimes }\widetilde{%
f}^{\varphi })(c\otimes 1)=\sum \widetilde{f}(1_{\psi })\widetilde{g}%
(c^{\psi })=\widetilde{f}(1_{A})\widetilde{g}(c)=(\varepsilon _{\widetilde{C}%
}(\widetilde{f})\widetilde{g})(c).
\end{equation*}
Hence $(\widetilde{A},\widetilde{C},\varphi )$ is a right-right entwining
structure. Since $(A,\widetilde{C})$ is a measuring $\alpha $-pairings, it
follows by Theorem \ref{equal} that $_{R}\widetilde{C}$ is locally
projective. The isomorphisms of categories \ref{du-ent-is} and \ref
{du-ent-is-fg} follow then by Theorem \ref{ent-sg}.$\blacksquare $
\end{Beweis}

\begin{lemma}
Consider the entwining structure $(\widetilde{A},\widetilde{C},\varphi ).$

\begin{enumerate}
\item  Consider the measuring $\alpha $-pairing $(A,\widetilde{C}).$ Let $%
M\in \mathcal{M}_{A}^{C}(\psi )$ and consider $M^{\ast }$ with the induced
right $\widetilde{A}$-module and left $A$-module structures. Then $M_{r}:=%
\mathrm{Rat}^{\widetilde{C}}(_{A}M^{\ast })\in \mathcal{M}_{\widetilde{A}}^{%
\widetilde{C}}(\varphi ).$

If $M,N\in \mathcal{M}_{A}^{C}(\psi )$ and $f:M\rightarrow N$ is $A$-linear $%
C$-colinear, then $f^{\ast }:N_{r}\rightarrow $ $M_{r}$ is $\widetilde{A}$%
-linear $\widetilde{C}$-colinear.

\item  Assume the measuring $R$-pairing $(\widetilde{A},C)$ to satisfy the $%
\alpha $-condition \emph{(}equivalently,\emph{\ }$_{R}C$ is locally
projective and $\widetilde{A}\subseteq C^{\ast }$ is dense\emph{)}. Let $%
K\in \mathcal{M}_{\widetilde{A}}^{\widetilde{C}}(\varphi )$ and consider $%
K^{\ast }$ with the induced left $\widetilde{A}$-module and right $A$-module
structures. Then $K^{r}:=\mathrm{Rat}^{C}(_{\widetilde{A}}K^{\ast })\in 
\mathcal{M}_{A}^{C}(\psi ).$

If $K,L\in \mathcal{M}_{\widetilde{A}}^{\widetilde{C}}(\varphi )$ and $%
g:K\rightarrow L$ is $\widetilde{A}$-linear $\widetilde{C}$-colinear, then $%
g^{\ast }:L^{r}\rightarrow $ $K^{r}$ is $A$-linear $C$-colinear.
\end{enumerate}
\end{lemma}

\begin{Beweis}
\begin{enumerate}
\item  Let $M\in \mathcal{M}_{A}^{C}(\psi ).$ Since $(A,\widetilde{C})$ is a
measuring $\alpha $-pairing, $M_{r}:=\mathrm{Rat}^{\widetilde{C}%
}(_{A}M^{\ast })$ is by Theorem \ref{equal} a right $\widetilde{C}$%
-comodule. Moreover we have for all $a\in A,$ $\widetilde{g}\in \widetilde{A}%
,$ $m\in M$ and $h\in M_{r}:$ 
\begin{equation*}
\begin{tabular}{lll}
$\lbrack a(h\widetilde{g})](m)$ & $=$ & $(h\widetilde{g})(ma)$ \\ 
& $=$ & $h(\widetilde{g}[ma])$ \\ 
& $=$ & $\sum h((ma)_{<0>}\widetilde{g}((ma)_{<1>}))$ \\ 
& $=$ & $\sum h(m_{<0>}a_{\psi }\widetilde{g}(m_{<1>}^{\psi }))$ \\ 
& $=$ & $\sum (a_{\psi }h)(m_{<0>})\widetilde{g}(m_{<1>}^{\psi })$ \\ 
& $=$ & $\sum h_{<0>}(m_{<0>})h_{<1>}(a_{\psi })\widetilde{g}(m_{<1>}^{\psi
})$ \\ 
& $=$ & $\sum h_{<0>}(m_{<0>})\widetilde{g}_{\varphi
}(m_{<1>})h_{<1>}^{\varphi }(a)$ \\ 
& $=$ & $\sum h_{<0>}(\widetilde{g}_{\varphi }m)h_{<1>}^{\varphi }(a)$ \\ 
& $=$ & $\sum (h_{<0>}\widetilde{g}_{\varphi })(m)h_{<1>}^{\varphi }(a)$ \\ 
& $=$ & $(\sum (h_{<0>}\widetilde{g}_{\varphi })h_{<1>}^{\varphi }(a))(m),$%
\end{tabular}
\end{equation*}
i.e. $h\widetilde{g}\in M_{r}$ with $\varrho (h\widetilde{g})=\sum h_{<0>}%
\widetilde{g}_{\varphi }\otimes h_{<1>}^{\varphi }.$ Hence $M_{r}\in 
\mathcal{M}_{\widetilde{A}}^{\widetilde{C}}(\varphi ).$

The second statement follows now by Lemma \ref{clos} (4) and Theorem \ref
{equal}.

\item  Let $K\in \mathcal{M}_{\widetilde{A}}^{\widetilde{C}}(\varphi ).$ By
assumption $(\widetilde{A},C)$ satisfies the $\alpha $-condition, hence $%
K^{r}:=\mathrm{Rat}^{C}(_{\widetilde{A}}K^{\ast })$ is by Theorem \ref{equal}
a right $C$-comodule. Moreover we have for all $a\in A,$ $\widetilde{g}\in 
\widetilde{A},$ $k\in K$ and $f\in K^{r}:$%
\begin{equation*}
\begin{tabular}{lll}
$\lbrack \widetilde{g}(fa)](k)$ & $=$ & $(fa)(k\widetilde{g})$ \\ 
& $=$ & $f(a(k\widetilde{g}))$ \\ 
& $=$ & $\sum f((k\widetilde{g})_{<0>}(k\widetilde{g})_{<1>}(a))$ \\ 
& $=$ & $\sum f(n_{<0>}\widetilde{g}_{\varphi }(k_{<1>}^{\varphi })(a))$ \\ 
& $=$ & $\sum (\widetilde{g}_{\varphi }f)(k_{<0>})k_{<1>}^{\varphi }(a)$ \\ 
& $=$ & $\sum f_{<0>}(k_{<0>})\widetilde{g}_{\varphi
}(f_{<1>})k_{<1>}^{\varphi }(a)$ \\ 
& $=$ & $\sum f_{<0>}(k_{<0>})<a_{\psi },k_{<1>}>\widetilde{g}(f_{<1>}^{\psi
})$ \\ 
& $=$ & $\sum f_{<0>}(a_{\psi }k)\widetilde{g}(f_{<1>}^{\psi })$ \\ 
& $=$ & $\sum (f_{<0>}a_{\psi })(k)\widetilde{g}(f_{<1>}^{\psi })$ \\ 
& $=$ & $(\sum (f_{<0>}a_{\psi })\widetilde{g}(f_{<1>}^{\psi }))(k),$%
\end{tabular}
\end{equation*}
i.e. $fa\in \mathrm{Rat}^{C}(_{\widetilde{A}}K^{\ast })$ with $\varrho
(fa)=\sum f_{<0>}a_{\psi }\otimes f_{<1>}^{\psi }.$ Hence $K^{r}\in \mathcal{%
M}_{A}^{C}(\psi ).$ As in (1), the second statement follows by Lemma \ref
{clos} (4) and Theorem \ref{equal}.$\blacksquare $
\end{enumerate}
\end{Beweis}

\begin{definition}
With the notation and assumptions above kept, we call the right-right
entwining structure $(\widetilde{A},\widetilde{C},\varphi )$ a \emph{dual
entwining structure}\textbf{\ }of $(A,C,\psi ).$ We also call $M_{r}$ (resp. 
$K^{r}$) a \emph{dual entwined module} of $M$ (resp. of $K$).
\end{definition}

\begin{theorem}
\label{main}Assume that $(\widetilde{A},C)$ satisfies the $\alpha $%
-condition. Then we have right adjoint contravariant functors 
\begin{equation*}
(-)_{r}:\mathcal{M}_{A}^{C}(\psi )\rightarrow \mathcal{M}_{\widetilde{A}}^{%
\widetilde{C}}(\varphi )\text{ and }(-)^{r}:\mathcal{M}_{\widetilde{A}}^{%
\widetilde{C}}(\varphi )\rightarrow \mathcal{M}_{A}^{C}(\psi ).
\end{equation*}
\end{theorem}

\begin{Beweis}
Let $M\in \mathcal{M}_{A}^{C}(\psi ),$ $K\in \mathcal{M}_{\widetilde{A}}^{%
\widetilde{C}}(\varphi )$ and consider the canonical $R$-linear maps 
\begin{equation*}
\lambda _{M}:M\rightarrow (M_{r})^{\ast }\text{ and }\lambda
_{K}:K\rightarrow (K^{r})^{\ast }.
\end{equation*}
Clearly $\lambda _{M}$ is $\widetilde{A}$-linear and $\lambda _{K}$ is $A$%
-linear, hence $\lambda _{M}(M)\subseteq (M_{r})^{r}$ and $\lambda
_{K}(K)\subseteq (K^{r})_{r}$ by Lemma \ref{clos} (4). It's easy then to see
that the right-adjointness is given by the functorial inverse isomorphisms 
\begin{equation*}
\begin{tabular}{lllllll}
$\Lambda _{M,K}:$ & $\mathrm{Hom}_{A}^{C}(M,K^{r})$ & $\rightarrow $ & $%
\mathrm{Hom}_{\widetilde{A}}^{\widetilde{C}}(K,M_{r}),$ & $f$ & $\mapsto $ & 
$\text{ }f^{\ast }\circ \lambda _{K},$ \\ 
$\Gamma _{M,K}:$ & $\mathrm{Hom}_{\widetilde{A}}^{\widetilde{C}}(K,M_{r})$ & 
$\rightarrow $ & $\mathrm{Hom}_{A}^{C}(M,K^{r}),$ & $g$ & $\mapsto $ & $%
g^{\ast }\circ \lambda _{M}.\blacksquare $%
\end{tabular}
\end{equation*}
\end{Beweis}

\begin{punto}
\label{dual-cat}Let $(A,C,\psi ),$ $(B,D,\Psi )$ be right-right entwining
structures and assume that $A,B$ are $\alpha $-algebras. Let $(\gamma
,\delta ):(A,C,\psi )\rightarrow (B,D,\Psi )$ be a morphism in $\mathbb{E}%
_{\bullet }^{\bullet },$ $\gamma ^{\circ }:B^{\circ }\rightarrow A^{\circ }$
the induced $R$-coalgebra morphism and $\delta ^{\ast }:D^{\ast }\rightarrow
C^{\ast }$\ the\ induced $R$-algebra morphism. Let $\widetilde{C}\subseteq
A^{\circ },$ $\widetilde{D}\subseteq B^{\circ }$ be $R$-subcoalgebras and $%
\widetilde{A}\subseteq C^{\ast },$ $\widetilde{B}\subseteq D^{\ast }$ be $R$%
-subalgebras with $\varepsilon _{C}\in \widetilde{A},$ $\varepsilon _{D}\in 
\widetilde{B}$ and assume that $\gamma ^{\circ }(\widetilde{D})\subseteq 
\widetilde{C}$ and $\delta ^{\ast }(\widetilde{B})\subseteq \widetilde{A}.$
Assume moreover that $\psi ^{\ast }(\widetilde{C}\otimes _{R}\widetilde{A}%
)\subseteq \widetilde{A}\otimes _{R}\widetilde{C},$ $\Psi ^{\ast }(%
\widetilde{D}\otimes _{R}\widetilde{B})\subseteq \widetilde{B}\otimes _{R}%
\widetilde{D}$ and let $(\widetilde{A},\widetilde{C},\varphi )$ and $(%
\widetilde{B},\widetilde{D},\Phi )$ be the induced dual entwining structures
of $(A,C,\psi )$ and $(B,D,\Psi )$ respectively. Then we have for all $%
\widetilde{g}\in \widetilde{B},$ $\widetilde{f}\in \widetilde{D},$ $d\in D$
and $b\in B:$%
\begin{equation*}
\begin{tabular}{lll}
$(\sum \delta ^{\ast }(\widetilde{g}_{\Phi })\underline{\otimes }\gamma
^{\circ }(\widetilde{f}^{\Phi }))(c\otimes a)$ & $=$ & $\sum (\widetilde{g}%
_{\Phi }\underline{\otimes }\widetilde{f}^{\Phi })(\delta (c)\otimes \gamma
(a))$ \\ 
& $=$ & $\sum (\widetilde{f}\underline{\otimes }\widetilde{g})(\gamma
(a)_{\Psi }\otimes \delta (c)^{\Psi })$ \\ 
& $=$ & $\sum (\widetilde{f}\underline{\otimes }\widetilde{g})(\gamma
(a_{\psi })\otimes \delta (c^{\psi }))$ \\ 
& $=$ & $\sum (\gamma ^{\circ }(\widetilde{f})\underline{\otimes }\delta
^{\ast }(\widetilde{g}))(a_{\psi }\otimes c^{\psi })$ \\ 
& $=$ & $\sum (\delta ^{\ast }(\widetilde{g})_{\varphi }\underline{\otimes }%
\gamma ^{\circ }(\widetilde{f})^{\varphi })(c\otimes a),$%
\end{tabular}
\end{equation*}
i.e. $(\delta ^{\ast },\gamma ^{\circ }):(\widetilde{B},\widetilde{D},\Phi
)\rightarrow (\widetilde{A},\widetilde{C},\varphi )$ is a morphism in $%
\mathbb{E}_{\bullet }^{\bullet }.$
\end{punto}

\section{Dual Doi-Koppinen modules}

\emph{Doi-Koppinen structures} were presented independently by Y. Doi \cite
{Doi92} and M. Koppinen \cite{Kop95} and provide a fundamental example of
entwining structures. The corresponding categories of Doi-Koppinen modules
unify themselves many categories of modules well studied by Hopf-algebraists
such as the categories of \emph{Hopf modules } \cite[4.1]{Swe69}, \emph{%
relative Hopf modules }\cite{Doi83},\emph{\ Doi's }$[C,H]$\emph{-modules} 
\cite{Doi83}, \emph{Dimodules,} \emph{Yetter-Drinfeld modules }and \emph{%
modules graded by }$G$\emph{-sets} \cite{CMZ2002}.

\section*{Dual module (co)algebras $\&$ comodule (co)algebras}

\qquad Before we present our dual Doi-Koppinen modules we introduce some
definitions and results concerning duality of (co)module (co)algebras.

\begin{definition}
Let $H$ be an $R$-bialgebra.

\begin{enumerate}
\item  A \emph{right }$H$\emph{-module algebra} is an $R$-algebra $(A,\mu
_{A},\eta _{A})$ with a right $H$-module structure through $\phi
_{A}:A\otimes _{R}H\rightarrow A,$ such that $\mu _{A}$ and $\eta _{A}$ are $%
H$\emph{-linear}, i.e. 
\begin{equation}
(ab)h=\sum (ah_{1})(bh_{2})\text{ and }1_{A}h=\varepsilon _{H}(h)1_{A}\text{
for all }a,b\in A\text{ and }h\in H.  \label{rma}
\end{equation}
Analog we define a \emph{left }$H$\emph{-module algebra}. An $H$\emph{%
-bimodule algebra,} is a left and a right $H$-module algebra, such that $A$
is an $H$-bimodule with the given left and right $H$-actions.

\item  A \emph{right }$H$\emph{-module coalgebra }is an $R$-coalgebra $%
(C,\Delta _{C},\varepsilon _{C})$ with a right $H$-module structure through $%
\phi _{C}:C\otimes _{R}H\rightarrow C,$ such that $\Delta _{C}$ and $%
\varepsilon _{C}$ are $H$\emph{-linear} (equivalently, $\phi _{C}$ is an $R$%
-coalgebra morphism), i.e. 
\begin{equation}
\Delta _{C}(ch)=\sum c_{1}h_{1}\otimes c_{2}h_{2}\text{ and }\varepsilon
_{C}(ch)=\varepsilon _{C}(c)\varepsilon _{H}(h)\text{ for all }c\in C\text{
and }h\in H.  \label{rmc}
\end{equation}
Analog we define a \emph{left }$H$\emph{-module coalgebra}. An $H$\emph{%
-bimodule coalgebra}, is a left and a right $H$-module coalgebra, which is
an $H$-bimodule with the given left and right $H$-actions.

\item  A \emph{right }$H$\emph{-comodule algebra} is an $R$-algebra $(A,\mu
_{A},\eta _{A})$ with a right $H$-comodule structure through $\varrho
_{A}:A\rightarrow A\otimes _{R}H$, such that $\mu _{A}$ and $\eta _{A}$ are $%
H$\emph{-colinear }(equivalently, $\varrho _{A}$ is an $R$-algebra
morphism), i.e. 
\begin{equation}
\varrho _{A}(ab)=\sum a_{<0>}b_{<0>}\otimes a_{<1>}b_{<1>}\text{ and }%
\varrho _{A}(1_{A})=1_{A}\otimes 1_{H}.  \label{rca}
\end{equation}
Analog we define a \emph{left }$H$\emph{-comodule algebra}. An $H$\emph{%
-bicomodule algebra} is an left and right $H$-comodule algebra, which is an $%
H$-bicomodule under the given left and right $H$-coactions.

\item  A \emph{right }$H$\emph{-comodule coalgebra} is an $R$-coalgebra $%
(C,\Delta _{C},\varepsilon _{C})$ with a right $H$-comodule structure
through $\varrho _{C}:C\rightarrow C\otimes _{R}H,$ such that $\Delta _{C}$
and $\varepsilon _{C}$ are $H$\emph{-colinear}, i.e. 
\begin{equation}
\sum c_{<0>1}\otimes c_{<0>2}\otimes c_{<1>}=\sum c_{1<0>}\otimes
c_{2<0>}\otimes c_{1<1>}c_{2<1>},\sum \varepsilon
_{C}(c_{<0>})c_{<1>}=\varepsilon _{C}(c)1_{H}.  \label{com-coal}
\end{equation}
Analog we define a \emph{left }$H$\emph{-comodule coalgebra}. An $H$\emph{%
-bicomodule coalgebra }is a left and a right $H$-comodule coalgebra, which
is an $H$-bicomodule with the given left and right $H$-coactions.
\end{enumerate}
\end{definition}

\begin{lemma}
\label{A0-module}Let $R$ be noetherian and $H$ an $R$-bialgebra. If $A$ is a
right \emph{(}resp. a left) $H$-module algebra, then $A^{\circ }\subset
A^{\ast }$ is a left \emph{(}resp. a right\emph{) }$H$-submodule. If $A$ is
an $H$-bimodule algebra, then $A^{\circ }\subseteq H^{\ast }$ is an $H$%
-subbimodule.
\end{lemma}

\begin{Beweis}
Let $A$ be a right $H$-module algebra. If $f\in A^{\circ },$ then we have
for all $h\in H$ and $a,b\in A:$%
\begin{equation*}
\begin{tabular}{lll}
$(b(hf))(a)$ & $=$ & $(hf)(ab)$ \\ 
& $=$ & $f((ab)h)$ \\ 
& $=$ & $f(\sum (ah_{1})(bh_{2}))$ \\ 
& $=$ & $\sum f_{1}(ah_{1})f_{2}(bh_{2})$ \\ 
& $=$ & $[\sum (h_{2}f_{2})(b)(h_{1}f_{1})](a).$%
\end{tabular}
\end{equation*}
So $hf\in A^{\circ }$ for every $h\in H,$ i.e. $A^{\circ }\subset A^{\ast }$
is a left $H$-submodule.

If $A$ is a left $H$-module algebra, then a similar argument shows that $%
A^{\circ }\subset A^{\ast }$ is a right $H$-submodule. The last statement is
obvious.$\blacksquare $
\end{Beweis}

\begin{proposition}
\label{cc-ma}Let $R$ be noetherian, $H$ an $\alpha $-bialgebra and $%
U\subseteq H^{\circ }$ an $R$-subbialgebra.

\begin{enumerate}
\item  Consider the measuring $R$-pairing $(U,H).$ If $A$ is a right \emph{(}%
a left\emph{)} $H$-comodule algebra, then $A$ is a left \emph{(}a right\emph{%
)}\ $U$-module algebra and $A^{\circ }$ is a right \emph{(}a left\emph{)} $U$%
-module coalgebra. If $A$ is an $H$-bicomodule algebra, then $A$ is a\ $U$%
-bimodule algebra and $A^{\circ }$ is a $U$-bimodule coalgebra.

\item  Consider the measuring $\alpha $-pairing $(H,U).$

\begin{enumerate}
\item  If $A$ is a right \emph{(}a left\emph{)} $U$-comodule algebra, then $%
A $ is a left \emph{(}a right\emph{) }$H$-module algebra. If $A$ is a $U$%
-bicomodule algebra, then $A$ is an $H$-bimodule algebra.

\item  If $A$ is a left \emph{(}a right\emph{) }$H$-module algebra, then $%
\mathrm{Rat}^{U}(_{H}A)$ \emph{(}resp. $^{U}\mathrm{Rat}(A_{H})$\emph{)} is
a right \emph{(}a left\emph{)} $U$-comodule algebra. If $A$ is an $H$%
-bimodule algebra, then $^{U}\mathrm{Rat}^{U}(_{H}A_{H})$ is a $U$%
-bicomodule algebra.
\end{enumerate}
\end{enumerate}
\end{proposition}

\begin{Beweis}
\begin{enumerate}
\item  W.l.o.g. assume that $A$ is a right $H$-comodule algebra through an $%
R $-algebra morphism $\varrho :A\rightarrow A\otimes _{R}H.$ For all $a,b\in
A$ and $f\in U$ we have 
\begin{equation*}
\begin{tabular}{lllll}
$f\rightharpoonup ab$ & $=$ & $\sum (ab)_{<0>}f((ab)_{<1>})$ & $=$ & $\sum
a_{<0>}b_{<0>}f(a_{<1>}b_{<1>})$ \\ 
&  & $\sum a_{<0>}b_{<0>}f_{1}(a_{<1>})f_{2}(b_{<1>})$ & $=$ & $\sum
(f_{1}\rightharpoonup a)(f_{2}\rightharpoonup b),$%
\end{tabular}
\end{equation*}
and moreover 
\begin{equation*}
f\rightharpoonup 1_{A}=\sum 1_{<0>}f(1_{<1>})=1_{A}f(1_{H})=1_{A}\varepsilon
_{U}(f).
\end{equation*}
Hence $A$ is a left $U$-module algebra.

Consider now the canonical $R$-linear map $\varpi :A^{\circ }\otimes
_{R}U\rightarrow (A\otimes _{R}H)^{\circ }.$ Then $(A,A^{\circ }),$ $%
(A\otimes _{R}H,A^{\circ }\otimes _{R}U)$ are measuring $\alpha $-parings
and we have a morphism of $R$-pairings 
\begin{equation*}
(\varrho ,\varrho ^{\circ }\circ \varpi ):(A\otimes _{R}H,A^{\circ }\otimes
_{R}U)\rightarrow (A,A^{\circ }).
\end{equation*}
Moreover $A^{\circ }\otimes _{R}A^{\circ }\hookrightarrow (A\otimes
_{R}A)^{\ast }$ and it follows from the assumption and Lemma \ref{prop-mes}
(1) that $\varrho ^{\circ }\circ \varpi :A^{\circ }\otimes _{R}U\rightarrow
A^{\circ }$ is an $R$-coalgebra morphism, i.e. $A^{\circ }$ is a right $U$%
-module coalgebra. If $A$ is an $H$-bicomodule, then $A$ is a\ $U$-bimodule
by Theorem \ref{equal} and $A^{\circ }\subseteq A^{\ast }$ is a $U$%
-subbimodule by Lemma \ref{A0-module}, hence $A$ is a $U$-bimodule algebra
and $A^{\circ }$ is a $U$-bimodule algebra.

\item  Consider the measuring $\alpha $-pairing $(H,U).$

\begin{enumerate}
\item  W.l.o.g. let $A$ be a right $U$-comodule algebra. Then we have for
all $h\in H$ and $a,b\in A:$%
\begin{equation*}
\begin{tabular}{lll}
$h\rightharpoonup (ab)$ & $=$ & $\sum (ab)_{<0>}<h,(ab)_{<1>}>$ \\ 
& $=$ & $\sum a_{<0>}b_{<0>}<h,a_{<1>}\star b_{<1>}>$ \\ 
& $=$ & $\sum a_{<0>}b_{<0>}<h_{1},a_{<1>}><h_{2},b_{<1>}>$ \\ 
& $=$ & $\sum (h_{1}\rightharpoonup a)(h_{2}\rightharpoonup b)$%
\end{tabular}
\end{equation*}
and 
\begin{equation*}
h\rightharpoonup 1_{A}=\varepsilon (h)1_{A},
\end{equation*}
\newline
i.e. $A$ is a left $H$-module algebra. If $A$ is a $U$-bicomodule algebra,
then $A$ is by Theorem \ref{equal} an $H$-bimodule, hence an $H$-bimodule
algebra.

\item  Assume now that $A$ is a left $H$-module algebra. Then we have for
all $a,b\in \mathrm{Rat}^{U}(_{H}A)$ and $h\in H:$%
\begin{equation*}
\begin{tabular}{lll}
$(h\rightharpoonup ab)$ & $=$ & $\sum (h_{1}\rightharpoonup
a)(h_{2}\rightharpoonup b)$ \\ 
& $=$ & $\sum a_{<0>}<h_{1},a_{<1>}>b_{<0>}<h_{2},b_{<1>}>$ \\ 
& $=$ & $\sum a_{<0>}b_{<0>}<h,a_{<1>}\star b_{<1>}>,$%
\end{tabular}
\end{equation*}
i.e. $ab\in \mathrm{Rat}^{U}(_{H}A)$ with $\varrho (ab)=\sum
a_{<0>}b_{<0>}\otimes a_{<1>}\star b_{<1>}.$ Note also that $%
h\rightharpoonup 1_{A}=\varepsilon _{H}(h)1_{A},$ i.e. $1_{A}\in \mathrm{Rat}%
^{U}(_{H}A),$ with $\varrho (1_{A})=1_{A}\otimes \varepsilon
_{H}=1_{A}\otimes 1_{U}.$ Hence $\mathrm{Rat}^{U}(_{H}A)$ is a right $U$%
-comodule algebra. If $A$ is an $H$-bimodule algebra, then $^{U}\mathrm{Rat}%
^{U}(_{H}A_{H})$ is by Theorem \ref{equal} a $U$-bicomodule, hence a $U$%
-bicomodule algebra.$\blacksquare $
\end{enumerate}
\end{enumerate}
\end{Beweis}

\begin{proposition}
\label{c-box}Let $R$ be noetherian, $H$ an $\alpha $-bialgebra and $%
U\subseteq H^{\circ }$ an $R$-subbialgebra.

\begin{enumerate}
\item  Consider the measuring $\alpha $-pairing $(H,U).$ If $C$ is a right 
\emph{(}a left\emph{)} $H$-module coalgebra, then $C^{\ast }$ is a left 
\emph{(}a right\emph{)}$\;H$-module algebra and $\mathrm{Rat}%
^{U}(_{H}C^{\ast })$ is a right \emph{(}a left\emph{)}\ $U$-comodule
algebra. If $C$ is an $H$-bimodule coalgebra, then $C^{\ast }$ is an $H$%
-bimodule algebra and $^{U}\mathrm{Rat}^{U}(_{H}C_{H}^{\ast })$ is a $U$%
-bicomodule algebra.

\item  Consider the measuring $R$-pairing $(U,H).$ If $C$ is a right \emph{(}%
a left\emph{)} $H$-comodule coalgebra, then $C$ is a left \emph{(}a right%
\emph{)} $U$-module coalgebra and $C^{\ast }$ is a right \emph{(}a left\emph{%
) }$U$-module algebra. If $C$ is an $H$-bicomodule coalgebra, then $C$ is a $%
U$-bimodule coalgebra and $C^{\ast }$ is a $U$-bimodule algebra.
\end{enumerate}
\end{proposition}

\begin{Beweis}
\begin{enumerate}
\item  Let $C$ be a right $H$-module coalgebra. Then we have for all $f,g\in
C^{\ast },$ $h\in H,$ $c\in C:$%
\begin{equation*}
\begin{tabular}{lllll}
$(h\rightharpoonup (f\star g))(c)$ & $=$ & $(f\star g)(ch)$ & $=$ & $\sum
f((ch)_{1})g((ch)_{2})$ \\ 
& $=$ & $\sum f(c_{1}h_{1})g(c_{2}h_{2})$ & $=$ & $\sum
(h_{1}f)(c_{1})(h_{2}g)(c_{2})$ \\ 
& $=$ & $(\sum (h_{1}f)\star (h_{2}g))(c)$ &  & 
\end{tabular}
\end{equation*}
and 
\begin{equation*}
(h\varepsilon _{C})(c)=\varepsilon _{C}(ch)=\varepsilon _{C}(c)\varepsilon
_{H}(h)=(\varepsilon _{H}(h)\varepsilon _{C})(c),
\end{equation*}
i.e. $C^{\ast }$ is a left $H$-module algebra. By Proposition \ref{cc-ma}
(2-b), $\mathrm{Rat}^{U}(_{H}C^{\ast })$ is a right $U$-comodule algebra. If 
$C$ is an $H$-bimodule coalgebra, then $C^{\ast }$ is an $H$-bimodule
algebra and $^{U}\mathrm{Rat}^{U}(_{H}C_{H}^{\ast })$ is a $U$-bicomodule by
Theorem \ref{equal}, hence a $U$-bicomodule algebra.

\item  W.l.o.g. assume that $C$ is a right $H$-comodule coalgebra. For all $%
c\in C,$ $f\in U$ we have 
\begin{equation*}
\begin{tabular}{lll}
$\sum (f\rightharpoonup c)_{1}\otimes (f\rightharpoonup c)_{2}$ & $=$ & $%
\sum (c_{<0>1}\otimes c_{<0>2}f(c_{<1>})$ \\ 
& $=$ & $\sum c_{1<0>}\otimes c_{2<0>}f(c_{1<1>}c_{2<1>})$ \\ 
& $=$ & $\sum c_{1<0>}\otimes c_{2<0>}f_{1}(c_{1<1>})f_{2}(c_{2<1>})$ \\ 
& $=$ & $\sum c_{1<0>}f_{1}(c_{1<1>})\otimes c_{2<0>}f_{2}(c_{2<1>})$ \\ 
& $=$ & $\sum f_{1}\rightharpoonup c_{1}\otimes f_{2}\rightharpoonup c_{2}$%
\end{tabular}
\end{equation*}
and 
\begin{equation*}
\varepsilon _{C}(f\rightharpoonup c)=\sum \varepsilon
_{C}(c_{<0>})f(c_{<1>})=f(\varepsilon _{C}(c)1_{H})=\varepsilon
_{C}(c)\varepsilon _{U}(f),
\end{equation*}
i.e. $C$ is a left $U$-module coalgebra. Analogous to (1)\ one can show that 
$C^{\ast }$ is a right $U$-module algebra. If $C$ is an $H$-bicomodule
coalgebra, then $C$ is a $U$-bimodule by Theorem \ref{equal}. Hence $C$ is a 
$U$-bimodule coalgebra and $C^{\ast }$ is a $U$-bimodule algebra.$%
\blacksquare $
\end{enumerate}
\end{Beweis}

\qquad The following result generalizes (\cite[Example 4.1.10.]{Mon93}):

\begin{corollary}
\label{dh-m} \smallskip Let $H\;$be an $R$-bialgebra and consider $H^{\ast }$
as an $H$-bimodule with the regular left and right $H$-actions.

\begin{enumerate}
\item  Since $H$ is an $H$-bimodule coalgebra, it follows \emph{(}by
Proposition \emph{\ref{c-box}(1))} that $H^{\ast }$ is an $H$-bimodule
algebra. If moreover $R$ is noetherian and $H$ is an $\alpha $-algebra, then 
$H^{\circ }\subset H^{\ast }$ is an $H$-subbimodule algebra.

\item  Let $R$ be noetherian, $H$ an $\alpha $-algebra and $U\subseteq
H^{\circ }$ an $R$-subbialgebra. Since $H$ is an $H$-bicomodule algebra, it
follows \emph{(}by Proposition \emph{\ref{cc-ma} (1))}, that $H$ is a $U$%
-bimodule algebra. In particular $H$ is an $H^{\circ }$-bimodule algebra.
\end{enumerate}
\end{corollary}

\section*{Doi-Koppinen modules}

\begin{punto}
\label{DK}A \emph{right-right Doi-Koppinen structure} over $R$ is a triple $%
(H,A,C)$ consisting of an $R$-bialgebra $H,$ a right $H$-comodule algebra $A$
and a right $H$-module coalgebra $C.$ Let $(H,A,C),$ $(K,B,D)$ be
right-right Doi-Koppinen structures. Then a morphism $(\beta ,\gamma ,\delta
):(H,A,C)\rightarrow (K,B,D)$ consists of an $R$-bialgebra morphism $\beta
:H\rightarrow K,$ an $R$-algebra morphism $\gamma :A\rightarrow B$ and an $R$%
-coalgebra morphism $\delta :C\rightarrow D,$ such that 
\begin{equation*}
\sum \gamma (a_{<0>})\otimes \delta (ca_{<1>})=\sum \gamma (a)_{<0>}\otimes
\delta (c)\gamma (a)_{<1>}\text{ for all }a\in A\text{ and }c\in C.
\end{equation*}
The category of right-right Doi-Koppinen modules is denoted by $\mathbb{DK}%
_{\bullet }^{\bullet }.$ For definitions of the categories of left-left,
right-left and left-right Doi-Koppinen structures the reader may refer to 
\cite{CMZ2002}.
\end{punto}

\begin{punto}
let $(A,H,C)$ be a right-right Doi-Koppinen structure. A \emph{right-right
Doi-Koppinen module} corresponding to $(H,A,C)$ is a right $A$-module $M,$
which is also a right $C$-comodule, such that 
\begin{equation*}
\varrho _{M}(ma)=\sum m_{<0>}a_{<0>}\otimes m_{<1>}a_{<1>}\text{ for all }%
m\in M\text{ and }a\in A.
\end{equation*}
For Doi-Koppinen modules $M,N$ corresponding to $(A,H,C)$ we denote with $%
\mathrm{Hom}_{A}^{C}(M,N)$ the set of all $A$-linear $C$-colinear maps from $%
M$ to $N.$ With $\mathcal{M}(H)_{A}^{C}$ we denote the category of
right-right Doi-Koppinen modules and $A$-linear $C$-colinear morphisms.
Setting 
\begin{equation}
\psi :C\otimes _{R}A\rightarrow A\otimes _{R}C,\text{ }c\otimes a\mapsto
\sum a_{<0>}\otimes ca_{<1>},  \label{psi-DK}
\end{equation}
it follows by \cite[Page 295]{Brz99} that $(A,C,\psi )$ is a right-right
entwining structure and $\mathcal{M}(H)_{A}^{C}\simeq \mathcal{M}%
_{A}^{C}(\psi ).$ Moreover $\#^{op}(C,A):=\mathrm{Hom}_{R}(C,A),$ introduced
first in \cite[2.2]{Kop95}, is an $R$-algebra with multiplication 
\begin{equation}
(f\cdot g)(c)=\sum f(c_{2})_{<0>}g(c_{1}f\left( c_{2}\right) _{<1>})
\label{DK-mult}
\end{equation}
and unity $\eta _{A}\circ \varepsilon _{C}.$
\end{punto}

\section*{The duality theorems}

In what follows we present for every $\alpha $-bialgebra $H$ over a
noetherian ground ring $R$ and every right $H$-module coalgebra (left $H$%
-module coalgebra) $C$ a right $H^{\circ }$-comodule $R$-algebra (a left $%
H^{\circ }$-comodule algebra) $C^{0},$ that plays an important role by the
dualization process in the rest of this note. In our \emph{infinite versions}
of duality theorems, $C^{0}$ will play the role of $C^{\ast }$ in the finite
versions (e.g. \cite{Yok82}).

\begin{definition}
Let $R$ be noetherian, $H$ an $\alpha $-bialgebra and consider the measuring 
$\alpha $-pairing $(H,H^{\circ }).$ For every right (resp. left)\ $H$-module
coalgebra $C$ we have by Proposition \ref{c-box} (1) the right (resp. the
left) $H^{\circ }$-comodule algebra 
\begin{equation*}
C^{0}:=\mathrm{Rat}^{H^{\circ }}(_{H}C^{\ast })\text{ (resp. }C^{0}:=\text{ }%
^{H^{\circ }}\mathrm{Rat}(C_{H}^{\ast })\text{).}
\end{equation*}
\qquad
\end{definition}

In view of \ref{du-co} and Propositions \ref{cc-ma}, \ref{c-box} we get

\begin{theorem}
\label{doi-hopf}Let $R$ be noetherian.

\begin{enumerate}
\item  Let $(H,A,C)$ be a right-right Doi-Koppinen structure and assume that 
$H,A$ are $\alpha $-algebras. Then $(H^{\circ },C^{0},A^{\circ })$ is a dual
right-right Doi-Koppinen structure of $(H,A,C)$ and we have isomorphism of
categories 
\begin{equation*}
\mathcal{M}(H^{\circ })_{C^{0}}^{A^{\circ }}\simeq \mathrm{Rat}^{A^{\circ }}(%
\mathcal{M}_{\#^{op}(A^{\circ },C^{0})})=\sigma \lbrack (C^{0}\otimes
_{R}A^{\circ })_{\#^{op}(A^{\circ },C^{0})}].
\end{equation*}
If moreover $_{R}A$ is f.g. projective, then 
\begin{equation*}
\mathcal{M}(H^{\circ })_{C^{0}}^{A^{\ast }}\simeq \mathcal{M}%
_{\#^{op}(A^{\ast },C^{0})}.
\end{equation*}

\item  Let $(\beta ,\gamma ,\delta ):(H,A,C)\rightarrow (K,B,D)$ be a
morphism in $\mathbb{DK}_{\bullet }^{\bullet }.$ If $H,K,A,B$ are $\alpha $%
-algebras and $\delta ^{\ast }(D^{0})\subseteq C^{0}$ \emph{(}e.g.\emph{\ }$%
\delta $ is $H$-linear, or $C^{0}=C^{\ast }$\emph{)}, then $(\beta ^{\circ
},\delta ^{0},\gamma ^{\circ }):(K^{\circ },D^{0},B^{\circ })\rightarrow
(H^{\circ },C^{0},A^{\circ })$ is a morphism in $\mathbb{DK}_{\bullet
}^{\bullet }.$
\end{enumerate}
\end{theorem}

\qquad As a corollary of Theorem \ref{main} we get the following

\begin{theorem}
\label{d-DK-m}Let $R$ be noetherian, $(H,A,C)$ a right-right Doi-Koppinen
structure with $H,A$ $\alpha $-algebras and consider the dual right-right
Doi-Koppinen structure $(H^{\circ },C^{0},A^{\circ }).$

\begin{enumerate}
\item  There is a contravariant functor 
\begin{equation*}
(-)^{0}:\mathcal{M}(H)_{A}^{C}\rightarrow \text{ }\mathcal{M}(H^{\circ
})_{C^{0}}^{A^{\circ }},\text{ }M\longmapsto M^{0}:=\mathrm{Rat}^{A^{\circ
}}(_{A}M^{\ast }).
\end{equation*}

\item  If $P:=(C^{0},C)$ satisfies the $\alpha $-condition \emph{(}%
equivalently $_{R}C$ is locally projective and $C^{0}\subset C^{\ast }$ is
dense), then there is a contravariant functor 
\begin{equation}
(-)^{\Diamond }:\mathcal{M}(H^{\circ })_{C^{0}}^{A^{\circ }}\rightarrow 
\mathcal{M}(H)_{A}^{C},\text{ }K\longmapsto K^{\Diamond }:=\mathrm{Rat}%
^{C}(_{C^{0}}K^{\ast }).  \label{box-DH}
\end{equation}
Moreover the contravariant functors $(-)^{0}$ and $(-)^{\Diamond }$ are
right adjoint.
\end{enumerate}
\end{theorem}

\begin{definition}
We say an $R$-algebra $A$ is $R$\emph{-c-cogenerated}, if for every $R$%
-cofinite ideal $I\vartriangleleft A,$ the $R$-module $A/I$ is $R$%
-cogenerated.
\end{definition}

\begin{remark}
Let $R$ be noetherian and $A$ an $\alpha $-algebra. For every right (resp.
left) $A$-module $M,$ set $M^{0}:=\mathrm{Rat}^{A^{\circ }}(_{A}M^{\ast })$
(resp. $M^{0}:=$ $^{A^{\circ }}\mathrm{Rat}(M_{A}^{\ast })$). If $A$ is $R$%
-c-cogenerated, then we have by \cite[Proposition 3.3.15]{Abu2001} 
\begin{equation}
M^{0}:=\{f\in M^{\ast }|\text{ }f(MI)=0\text{ (resp. }f(IM)=0\text{) for
some }R\text{-cofinite ideal }I\vartriangleleft A\}.  \label{C0}
\end{equation}
\end{remark}

\begin{example}
Let $A$ be an $R$-algebra, $C$ an $R$-coalgebra and consider the category of
right $A$-modules and right $C$-comodules satisfying the compatibility
relation 
\begin{equation*}
\varrho _{M}(ma)=\sum m_{<0>}a\otimes m_{<1>}\text{ for all }m\in M\text{
and }a\in A.
\end{equation*}
The category of such modules and $A$-linear $C$-colinear morphisms is called
the category of \emph{Long dimodules,} denoted by $\mathcal{L}_{A}^{C},$ and
was introduce by F. Long in \cite{Lon74}. Considering $A$ as a trivial $R$%
-comodule algebra and $C$ as a trivial right $R$-module coalgebra we get the
right-right Doi-Koppinen structure $(R,A,C)$ and it follows that $\mathcal{L}%
_{A}^{C}\simeq \mathcal{M}(R)_{A}^{C}.$ If $A$ is an $\alpha $-algebra, then 
$(R,C^{\ast },A^{\circ })$ is a dual right-right Doi-Koppinen structure of $%
(R,A,C)$ and the contravariant functors $(-)^{0}:\mathcal{L}%
_{A}^{C}\rightarrow \mathcal{L}_{C^{\ast }}^{A^{\circ }}$ and $(-)^{\Diamond
}:\mathcal{L}_{C^{\ast }}^{A^{\circ }}\rightarrow \mathcal{L}_{A}^{C}$ are
right adjoint.
\end{example}

Inspired by \cite{Liu94} and in contradiction to \cite[Page 138]{Abe80}, the
following example shows that for a Hopf $R$-algebra $H$ and an $H$-module
algebra $A$ over a field, the dual $R$-coalgebra $A^{\circ }$ need \emph{not}
be an $H^{\circ }$-comodule coalgebra.

\begin{c-example}
\label{geg-abe}Let $R$ be a field and $H$ a \emph{coreflexive} Hopf $R$%
-algebra with $\mathrm{dim}(H)=\infty $ (e.g. the Hopf $R$-algebra of 
\cite[Example 5]{Lin77}). By Lemma \ref{dh-m} $H^{\ast }$ is a right $H$%
-module algebra. If $H^{\ast \circ }\simeq H$ were a right $H^{\circ }$%
-comodule coalgebra, then we would have an $R$-cofinite ideal $%
J\vartriangleleft H$ with 
\begin{equation*}
0=<1_{H},H^{\ast }\leftharpoonup J>=<J,H^{\ast }>.
\end{equation*}
But we would get then $J=0$ (which contradicts the assumption $\mathrm{dim}%
(H)=\infty $).
\end{c-example}

\begin{remark}
Let $H$ be an $R$-bialgebra, $A$ a right $H$-module algebra and $C$ a right $%
H$-comodule coalgebra. Then $(H,A,C)$ is called a \emph{right-right
alternative Doi-Koppinen structure}. Such structures were introduced by P.
Schauenburg in \cite{Sch2000}, who showed that with 
\begin{equation*}
\psi :C\otimes _{R}A\rightarrow A\otimes _{R}C,\text{ }c\otimes a\mapsto
\sum ac_{<1>}\otimes c_{<0>},
\end{equation*}
$(A,C,\psi )$ is a right-right entwining structure. Moreover he gave an
example of such a structure that \emph{can not }be derived form a
Doi-Koppinen structure. The previous counter example shows that, even over
base fields, $(H^{\circ },C^{0},A^{\circ })$ \emph{may not} be a dual
alternative Doi-Koppinen structure of $(H,A,C).$
\end{remark}

\section*{Cleft $H$-Extensions}

\qquad Hopf-Galois extensions were presented by S. Chase and M. Sweedler 
\cite{CS69} for a \emph{commutative }$R$-algebra acting on a Hopf $R$-Hopf
and are considered as generalization of the classical Galois extensions over
fields (e.g. \cite[8.1.2]{Mon93}). In \cite{KT81} H. Kreimer and M. Takeuchi
extended these to the \emph{noncommutative }case.

\begin{punto}
$H$\textbf{-Extensions. }(\cite{Doi85}) Let $H$ be an $R$-bialgebra, $B$ a
right $H$-comodule algebra and $A:=B^{coH}=\{a\in B\mid \varrho (a)=a\otimes
1_{H}\}.$ Then $A$ is an $R$-algebra and the algebra extension $%
A\hookrightarrow B$ is called a \emph{right }$H$\emph{-extension.}

A (\emph{total}) \emph{integral} for $B$ is an $H$\emph{-colinear} map $%
\gamma :H\rightarrow B$ (with $\gamma (1_{H})=1_{B}$). If $B$ admits an
integral, that is invertible in $(\mathrm{Hom}_{R}(H,B),\star ),$ then $%
A\hookrightarrow B$ is called a \emph{cleft}\textbf{\ }right $H$-extension.
\end{punto}

\begin{example}
Let $H$ be an $R$-bialgebra. By \cite[Corollary 6]{DT86} $H/R$ is a \emph{%
cleft} $H$-extension, iff $H$ is a Hopf $R$-algebra. In this case $%
id_{H}:H\rightarrow H$ is an invertible total integral with inverse the
antipode $S_{H}.$
\end{example}

\begin{punto}
$H$-\textbf{Coextensions}. Let $H$ be an $R$-bialgebra and $D\;$a right $H$%
-module coalgebra. Then $H^{+}:=\mathrm{Ke}(\varepsilon _{H})$ in an $H$%
-coideal, $DH^{+}$ is a $D$-coideal and $C:=D/DH^{+}$ is a right $H$-module
coalgebra with the induced $H$-module structure. The canonical coalgebra
epimorphism $\pi :D\rightarrow C$ is called a right $H$\emph{-coextension}
of $D.$

A (\emph{total}) \emph{cointegral} for $D$ is an $H$-linear map $\omega
:D\rightarrow H$ (with $\varepsilon _{H}\circ \omega =\varepsilon _{D}$). A
right $H$-coextension $\pi :D\rightarrow C$ is called \emph{cocleft}, if $D$
admits cointegral, that is invertible in $(\mathrm{Hom}_{R}(D,H),\star ).$
\end{punto}

\qquad As a corollary of our results in this section we get

\begin{proposition}
\label{co-cleft}Let $R$ be noetherian, $H$ a Hopf $\alpha $-algebra with
bijective antipode, $D$ a right $H$-module coalgebra and $C:=D/DH^{+}.$ If $%
\pi :D\rightarrow C$ is a \emph{(cocleft)} right{\normalsize \ }$H$%
-coextension, then $\pi ^{\circ }:C^{0}\hookrightarrow D^{0}$ is a \emph{%
(cleft)} right $H^{\circ }$-extension.
\end{proposition}

\begin{Beweis}
Let $D$ be a right $H$-module coalgebra through $\phi _{D}:D\otimes
_{R}H\rightarrow D.$ By Proposition \ref{c-box} (1) $D^{0}$ is a right $%
H^{\circ }$-comodule algebra through $\phi _{D}^{\circ }:D^{0}\rightarrow
D^{0}\otimes _{R}H^{\circ }.$ Moreover we have $C^{\ast }=(D^{\ast
})^{H}:=\{g\in D^{\ast }\mid $ $hg=\varepsilon (h)g\mathbb{\ }$for every $%
h\in H\}.$ Hence $(D^{0})^{coH^{\circ }}=(D^{0})^{H}=$ $C^{0}$ (by 
\cite[Lemma 2.5.15]{Abu2001}), i.e. $\pi ^{\circ }:C^{0}\hookrightarrow
D^{0} $ is a right $H^{\circ }$-extension.

If $\omega :D\rightarrow H$ is a cointegral for $D,$ then $\omega $ is by
definition $H$-linear and so $\omega ^{\circ }\in \mathrm{Hom}_{H-}(H^{\circ
},D^{0})=\mathrm{Hom}^{H^{\circ }}(H^{\circ },D^{0}),$ i.e. $\omega ^{\circ
} $ is an integral for $D^{0}.$

Let $\omega $ be invertible in $(\mathrm{Hom}_{R}(D,H),\star )$ with inverse 
$\omega ^{-1}:D\rightarrow H.$ Analog to \cite{Zha98} we get 
\begin{equation*}
\omega ^{-1}(dh)=S_{H}(h)\omega ^{-1}(d)\text{ for all }h\in H\text{ and }%
d\in D.
\end{equation*}
If $f\in H^{\circ },$ then we have for all $d\in D$ and $h\in H:$%
\begin{equation*}
\begin{tabular}{lllll}
$(h(\omega ^{-1})^{\circ }(f))(d)$ & $=$ & $((\omega ^{-1})^{\circ }(f))(dh)$
& $=$ & $f(\omega ^{-1}(dh))$ \\ 
& $=$ & $f(S(h)\omega ^{-1}(d))$ & $=$ & $\sum f_{1}(S(h))f_{2}(\omega
^{-1}(d))$ \\ 
& $=$ & $\sum (S^{\circ }(f_{1}))(h)((\omega ^{-1})^{\circ }(f_{2}))(d)$ & $%
= $ & $(\sum S^{\circ }(f_{1})(h)(\omega ^{-1})^{\circ }(f_{2}))(d),$%
\end{tabular}
\end{equation*}
i.e. $(\omega ^{-1})^{\circ }\in D^{0}$ with $\varrho ((\omega ^{-1})^{\circ
})=\sum (\omega ^{-1})^{\circ }(f_{2})\otimes S^{\circ }(f_{1})$. Moreover
we have for all $f\in H^{\circ }$ and $d\in D:$%
\begin{equation*}
\begin{tabular}{lllll}
$((\omega ^{\circ }\star (\omega ^{-1})^{\circ })(f))(d)$ & $=$ & $((\omega
^{\circ }\underline{\otimes }(\omega ^{-1})^{\circ })(\Delta (f)))(d)$ &  & 
\\ 
& $=$ & $(\sum \omega ^{\circ }(f_{1})\star (\omega ^{-1})^{\circ
}(f_{2}))(d)$ &  &  \\ 
& $=$ & $(\sum \omega ^{\circ }(f_{1})\underline{\otimes }(\omega
^{-1})^{\circ }(f_{2}))(d_{1}\otimes d_{2})$ &  &  \\ 
& $=$ & $\sum \omega ^{\circ }(f_{1})(d_{1})(\omega ^{-1})^{\circ
}(f_{2})(d_{2})$ &  &  \\ 
& $=$ & $\sum f_{1}(\omega (d_{1}))f_{2}(\omega ^{-1}(d_{2}))$ &  &  \\ 
& $=$ & $\sum f(\omega (d_{1})\omega ^{-1}(d_{2}))$ &  &  \\ 
& $=$ & $f((\omega \star \omega ^{-1})(d))$ &  &  \\ 
& $=$ & $f(\varepsilon _{D}(d)1_{H})$ &  &  \\ 
& $=$ & $\varepsilon _{H^{\circ }}(f)\varepsilon _{D}(d),$ &  & 
\end{tabular}
\end{equation*}
i.e. $\omega ^{\circ }\star (\omega ^{-1})^{\circ }=id_{\mathrm{Hom}%
_{R}(H^{\circ },D^{0})}.$ Analog one can prove that $(\omega ^{-1})^{\circ
}\star \omega ^{\circ }=id_{\mathrm{Hom}_{R}(H^{\circ },D^{0})}.$ So $\omega
^{\circ }$ is $\star $-invertible with inverse $(\omega ^{-1})^{\circ }$ and 
$\pi ^{\circ }:C^{0}\hookrightarrow D^{0}$ is a \emph{cleft }right $H^{\circ
}$-extension.$\blacksquare $
\end{Beweis}

\end{document}